%% file: bun
\documentclass[11pt]{amsart}

\usepackage{amssymb, amsmath, epsfig, graphics, amscd, ifthen}

\newcommand{\pics}{true}

\newcommand{\susy}{{$N=1$}}
\newtheorem{theorem}{Theorem}[section]

\newtheorem{proposition}[theorem]{Proposition}
\newtheorem{lemma}[theorem]{Lemma}
\newenvironment{pf}{\paragraph{\sc Proof}}{\hspace*{\fill}$\Box$\par\medskip}

\newtheorem{remark}[theorem]{Remark}
\newtheorem{definition}[theorem]{Definition}

\newcommand{\spec}{\mathop{\rm Spec\,}\nolimits}
\newcommand{\coker}{\mathop{\rm coker}\nolimits}
\newcommand{\coh}{\mathop{\rm Coh\,}\nolimits}
\renewcommand{\mod}{\mathop{\rm Mod\,}\nolimits}
\newcommand{\proj}{\mathop{\rm Proj\,}\nolimits}

\newcommand{\Coh}{\mathop{\rm Coh}\nolimits}

\renewcommand{\phi}{\varphi}

\newcommand{\Hom}{\mathop{\rm Hom}\nolimits}

\newcommand{\Ext}{\mathop{\rm Ext}\nolimits}
\newcommand{\shExt}{\mathop{{\mathcal E}\hspace{-1pt}xt}\nolimits}

\newcommand{\SL}{{\rm SL}}
\newcommand{\GL}{{\rm GL}}
\newcommand{\SU}{{\rm SU}}

\renewcommand{\Box}{\square}

\newcommand{\la}{\lambda}

\newcommand{\F}{\mathcal F}

\newcommand{\M}{\mathcal M}
\newcommand{\B}{\mathcal B}
\newcommand{\D}{\mathcal D}

\newcommand{\I}{\mathcal I}
\newcommand{\J}{\mathcal J}
\newcommand{\A}{{\bf A}}
\renewcommand{\AA}{\mathbb A}
\newcommand{\sA}{\mathcal A}
\newcommand{\sP}{\mathcal P}
\newcommand{\X}{\mathcal X}
\newcommand{\T}{\mathcal T}
\newcommand{\Y}{\mathcal Y}
\newcommand{\E}{\mathcal E}
\newcommand{\R}{\mathbb R}
\newcommand{\U}{\mathcal U}
\newcommand{\V}{\mathcal V}
\newcommand{\W}{\mathcal W}
\newcommand{\C}{\mathbb C}
\newcommand{\sC}{\mathcal C}
\newcommand{\Z}{\mathbb Z}
\renewcommand{\O}{\mathcal O}

\newcommand{\PP}{\mathbb P}
\renewcommand{\P}{{\bf P}}

\newcommand{\sQ}{\mathcal Q}

\newcommand{\sH}{{\mathcal H}} 
\newcommand{\sh}{{\mathcal H}_0} 
\newcommand{\lH}{{\mathfrak h}} 
\newcommand{\lh}{{\mathfrak h}_0} 
\newcommand{\De}{\Delta} 
\newcommand{\de}{\Delta_0} 
\newcommand{\roo}{\eta} 
\newcommand{\posro}{R_+} 
\newcommand{\spvec}{{\bf v_0}} 

\newcommand{\defpar}{s}
\newcommand{\dpar}{s}
\newcommand{\defpoly}{\Theta}
\newcommand{\defspace}{{\mathcal S}}

\newcommand{\color}[6]{}

\begin{document}
\thispagestyle{empty}
\begin{center}
{\LARGE Sheaves on fibered threefolds and quiver sheaves} 

\vspace{0.2in}

{\large Bal\'azs Szendr\H oi} 

\vspace{0.15in}

{\large August 2007}

\end{center}

\vspace{0.1in}

{\small
\begin{center} {\sc Abstract} \end{center}
{\leftskip=50pt \rightskip=50pt
\noindent This paper classifies a class of holomorphic $D$-branes, closely related 
to framed torsion-free sheaves, on threefolds fibered in resolved $ADE$ surfaces 
over a general curve~$C$, in terms of representations with relations of a twisted 
Kronheimer--Nakajima-type quiver in the category $\Coh(C)$ of coherent sheaves 
on~$C$. For the local Calabi--Yau case $C\cong\A^1$ and special choice of 
framing, one recovers the \susy\ $ADE$ quiver studied by Cachazo--Katz--Vafa. 
\par}}

\section*{Introduction}

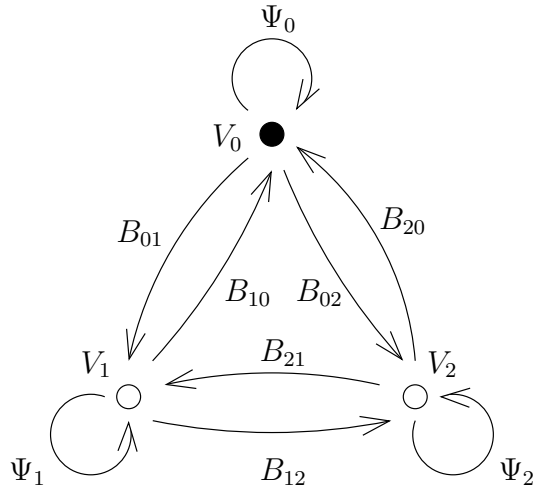
\begin{figure}[ht]
\centering
\ifthenelse{\equal{\pics}{true}}{\input{quiver2.pstex_t}}{}
\caption{A representation of the affine \susy\ $A_2$ quiver} 
\label{fig!quiv2}
\end{figure}

The purpose of this paper is to study, via dimensional reduction, certain 
holomorphic $D$-branes, closely related to torsion-free sheaves, on 
threefolds $X\to C$ fibered in resolved $ADE$ surfaces over a curve. 
Fibered local Calabi--Yau threefolds $X\rightarrow \AA^1$ of this type, 
as well as their deformations $X_\defpar\rightarrow \AA^1$ and extremal transitions,
were thoroughly analized in~\cite{ckv, cifkv} from the point of view of 
supersymmetric gauge theory. The paper~\cite{ckv} contains an assertion, made 
explicit in~\cite{katz} and studied in~\cite{zhu}, that exceptional components
of a natural threefold contraction $X_\defpar\to \bar X_\defpar$ are classified 
by irreducible representations of a certain quiver with loop edges, 
the \susy\ $ADE$ quiver (see Figure~\ref{fig!quiv2} for an example), 
satisfying a specific set of relations. This statement is in the spirit of 
Gabriel's theorem classifying exceptional (not necessarily irreducible) 
rational curves in resolved $ADE$ surfaces in terms of irreducible
representations of the corresponding Dynkin quiver.

In this paper we generalize the work of~\cite{ckv, katz, zhu} 
in two directions: we consider holomorphic $D$-branes, 
objects in the derived category of coherent sheaves,
instead of exceptional components, and we study the semi-local case: 
the neighbourhood of a deformed $ADE$ fibration $X_\defpar\rightarrow C$ over 
a general curve~$C$. The main result is Theorem~\ref{thm!class}, which shows that 
certain holomorphic $D$-branes on the fibered threefold $X_\defpar$ are classified 
by representations with relations of a Kronheimer--Nakajima-type quiver 
in the category $\Coh(C)$ of coherent sheaves on the curve~$C$. In particular, 
moduli spaces of such holomorphic $D$-branes are quiver bundle varieties over $C$.
If $C\cong\AA^1$, a further dimensional reduction leads to Theorem~\ref{thm!A1}, 
relating sheaves on the threefold to the zero-dimensional problem of ordinary 
matrix representations of the \susy\ $ADE$ quiver of~\cite{ckv, katz, zhu}.
The loops in the \susy\ $ADE$ quiver arise as the action by multiplication 
of a parameter $t\in H^0(\O_{\AA^1})$ on spaces of sections of sheaves on 
the base~$\AA^1$.

The geometry considered in this paper is non-monodromic, meaning that there
is no global~\cite{sz} nor local~\cite{ckv} monodromy in the fibration of
$ADE$ surfaces over the curve~$C$. It appears to be an interesting question
to extend the results proved here to these more general cases involving 
monodromy.

In recent work~\cite{penn}, the moduli space of certain very special 
holomorphic $D$-branes on resolved $A_1$-fibered geometries $X\rightarrow C$ 
has been connected, via imposing a superpotential and going through a large~$N$ 
transition, to the Hitchin system on~$C$. The branes studied in~\cite{penn} are 
not of the type classified by our results; they should rather correspond to a 
complex of quiver representations. Understanding the precise connection 
between~\cite{penn} and the present paper is left for future work.

After introducing basic notation in Section~\ref{sec!surf}, 
Section~\ref{sec!defs} describes the threefolds we study, and defines
some auxiliary sheaves of non-commutative algebras over the curve~$C$. 
Section~\ref{sec!res} contains our results, in particular the
general statement Theorem~\ref{thm!class} connecting quiver bundles 
to holomorphic $D$-branes on ADE fibrations, as well as the 
statement for the affine case. Proofs are discussed in Section~\ref{sec!pfs}. 

\section{Finite groups of type $ADE$ and surfaces}
\label{sec!surf}

Let $\Gamma<\SL(2,\C)$ be a finite subgroup of type~$A, D$ or~$E$. 
Let~$\lh$ be the Cartan subalgebra of the finite dimensional 
Lie algebra of the same type. Fix a set of simple roots $\{\roo_a: a\in\de\}$
indexed by nodes of the Dynkin diagram~$\de$, and let~$\posro$ be 
the set of positive roots. Let $\lH$ be the corresponding affine 
Cartan with simple roots indexed by nodes of the Dynkin diagram $\De\supset \de$. 

The group ring $\C\Gamma$ has center $Z(\C\Gamma)\cong\C^\De$; 
explicitly, for $\lambda \in Z(\C\Gamma)$, the isomorphism is obtained by 
taking the trace of $\lambda$ on a set of 
irreps, indexed by the nodes of $\De$ according to the McKay correspondence.   
There is also a natural identification
\[\lh=\{\lambda\in\C^\De\, |\, \lambda\cdot\delta=0\}\subset \lH\cong\C^\De,\]
where $\delta=(\delta_a)$ 
are the dimensions of the irreps of $\Gamma$.

\begin{lemma} The centralizer $C_{\GL(2,\C)}(\Gamma)$ of $\Gamma$ in 
${\GL(2,\C)}$ is 
\begin{enumerate}
\item the full group $\GL(2,\C)$ for type $A_1$;
\item a torus $(\C^\star)^2$ in~$\GL(2,\C)$ for type $A_n$ with $n>1$;
\item the center $\C^\star$ of~$\GL(2,\C)$ for types $D$ and $E$.
\end{enumerate}
\label{lemma!cent}
\end{lemma}

Let $\bar Y=\AA^2/\Gamma$ be the singular affine quotient,  
$Y\rightarrow \bar Y$ its minimal resolution. Exceptional curves in the 
resolution are in one-to-one correspondence with the nodes of $\de$,
and thus with a set of simple roots of $\lh$; 
the positive roots $\roo\in\posro$ correspond to connected, 
possibly reducible exceptional rational curves. The 
universal deformations $\Y\rightarrow\lh$ and $\bar\Y\rightarrow\lh/W$ of 
$Y$ and $\bar Y$, where $W$ denotes the Weyl group, 
are connected by the well known commutative diagram
\[\begin{array}{ccccc}
\Y & \longrightarrow & p^*\bar\Y & \longrightarrow &\bar\Y \\
&\searrow & \downarrow&&\downarrow\\
&& \lh & \stackrel{p}\longrightarrow & \lh/W.
\end{array}
\]

\section{Threefolds: definitions} \label{sec!defs} 
\subsection{The geometry}
Let $C$ be a curve, and let $\sQ$ be a rank-two
vector bundle on $C$ whose structure group reduces from $\GL(2,\C)$ to the
centralizer $C_{\GL(2,\C)}(\Gamma)$. Thus, by Lemma~\ref{lemma!cent}, 
\begin{itemize}
\item for type $A_1$, $\sQ$ is an arbitrary rank-two vector bundle; 
\item for type $A_n$ with $n>1$, $\sQ\cong\sQ_1\oplus\sQ_2$ 
is the direct sum of two line bundles;
\item for types $D,E$, $\sQ\cong\sQ_0^{\oplus 2}$ for some line bundle $\sQ_0$.
\end{itemize}
There is a fiberwise $\Gamma$-action on the total space of the 
vector bundle $\sQ$, and the quotient $\bar{X}=\sQ/\Gamma$ is a threefold 
with a curve of compound Du Val singularities
along the image of the zero section. 
Let $f\colon X\rightarrow \bar X$ be the crepant resolution,
with a map $\pi\colon X\rightarrow C$ whose fibres are minimal resolutions 
of the corresponding surface singularity, with trivial monodromy in the fibres. 
The canonical bundle of $X$ is
\[\omega_X\cong \pi^*(\omega_C\otimes \det\sQ^\vee).\]
In particular,~$X$ is Calabi--Yau if and only if~$\sQ$ has canonical 
determinant on~$C$. 

Part of the deformation theory of the threefold~$X$ was described 
in~\cite{sz}. Let $\sh=\det\sQ\otimes\lh$, a vector bundle over $C$, and let
$\defspace=H^0(C, \sh)$ be its space of sections. Then there is a smooth family
of threefolds $\X\rightarrow \defspace$, with injective Kodaira--Spencer 
map and central fibre $X_0\cong X$, together with a fibration 
$\X\to C\times S$ and a contraction $\X\to\bar\X$ over $S$. 
Thus, for every $\defpar\in\defspace$, the threefold fibre $X_\defpar$ possesses 
a fibration $\pi_\dpar\colon X_\defpar\rightarrow C$ in surfaces and a contraction 
$f_\defpar\colon X_\defpar\rightarrow \bar X_\defpar$ to a singular threefold 
with compound Du Val singularities. 
More precisely, for every positive root $\roo\in\posro$ of $\lh$, 
there is a map $p_\roo\colon \sh\rightarrow \det\sQ$, whose vanishing locus 
is a family of root hyperplanes in the $\lh$ fibers, and we have 

\begin{lemma} Let $\defpar\in \defspace=H^0(C,\sh)$ be a section of $\sh$, and 
let $\roo\in\posro$ be a positive root of $\lh$. 
The contraction $f_\defpar\colon X_\defpar\rightarrow \bar X_\defpar$ 
contracts a (possibly reducible) rational curve corresponding to the 
root~$\roo$ over a point~$P\in C$, if and only if the projected section
$p_\roo(\defpar)\in H^0(C, \det\sQ)$ vanishes at~$P\in C$. 
\end{lemma} 
\noindent
Thus if the projected section $p_\roo(\defpar)$ is not identically zero
for any root~$\roo$, 
then $f_\defpar$ is a small contraction, contracting rational curves to 
isolated singularities in certain configurations. If for different 
roots~$\roo$, the sections~$p_\roo(\defpar)$ have different simple zeros, 
then $f_\defpar$ contracts a set of isolated $(-1, -1)$-curves
to simple nodes. If the linear system~$\det\sQ$ has no base points on $C$,
then this holds for generic~$\defpar\in\defspace$.

In the special case~$C\cong\AA^1$, the central fiber~$X_0=\AA^1\times Y$ is 
Calabi--Yau, and its deformations are parameterized by an $\lh$-valued 
polinomial $\defpar\in\lh[t]$. Under the isomorphism
$\lh\cong\{\la\, | \, \defpar\cdot \delta=0\}\subset\C^\De$, we can also
parameterize deformations by a set of ordinary polinomials 
$\defpoly_a\in \C[t]$ 
indexed by nodes of the affine Dynkin diagram $\De$, satisfying 
$\sum_a \delta_a \defpoly_a=0$. The exceptional fibres of 
$f_\defpar\colon X_\defpar\rightarrow \bar X_\defpar$
lie over roots of the various polynomials $\defpoly_{\roo_a}=\defpoly_a$, 
corresponding to simple roots $\roo_a$, as well as over roots of their linear 
combinations $\defpoly_\roo=\sum_a \mu_a \defpoly_a$, corresponding to other 
positive roots $\roo=\sum_a\mu_a\roo_a\in\posro$. 
For generic choice of parameter $\defpar\in\defspace$, equivalently for 
generic choice of $\{\defpoly_a\}$, 
the polynomials $\{\defpoly_\roo: \roo\in\posro\}$ have 
distinct simple roots, and the exceptional set of 
$f_\defpar\colon X_\defpar\rightarrow \bar X_\defpar$ 
consists of isolated $(-1, -1)$-curves.

\subsection{Sheaves of non-commutative algebras and their sheaves of modules} 

Given $(C, \sQ)$, let $\sH=\det\sQ\otimes \lH$, a vector bundle on the 
curve~$C$ containing~$\sh$ as a subbundle. 
Given a section $\defpar\in H^0(C, \sH)$, consider the natural composition
\[ \sigma_\defpar\colon\ \  \sQ^\vee\otimes\sQ^\vee\stackrel{\wedge^2}\longrightarrow \det\sQ^\vee \stackrel{\cdot \defpar}\longrightarrow \lH\otimes \O_C\stackrel{\sim}\longrightarrow Z(\C\Gamma)\otimes\O_C,\]
a family of $Z(\C\Gamma)$-valued symplectic forms in the fibres of the vector 
bundle $\sQ^\vee$. Also fix, once and for all, a trivializing section 
$z\in H^0(\O_C)$. 

\begin{definition}\rm Let $\sA_\defpar$ be the sheaf of non-commutative algebras 
on $C$ whose sections on an open set $U\subset C$ are 
\[\sA_\defpar(U)= T\sQ^\vee(U)*\C\Gamma \big/\big\langle\!\big\langle[x_1,x_2]+\sigma_\defpar(x_1,x_2)\big\rangle\!\big\rangle,\]
where $T\sQ^\vee(U)$ is the full tensor algebra of $\sQ^\vee(U)$, 
$x_i\in\sQ^\vee(U)$ are local sections, and 
$\big\langle\!\big\langle\ldots\big\rangle\!\big\rangle$ denotes the two-sided
ideal generated by all given expressions. Define also
\[\sP_\defpar(U)=T(\sQ^\vee\oplus\O_C)(U)*\C\Gamma \big/\big\langle\!\big\langle[x_1,x_2]+\sigma_\defpar(x_1,x_2)z^2, [x_i, z]\big\rangle\!\big\rangle,\]
where the fixed section $z\in H^0(\O_C)$ commutes with elements 
of~$\C\Gamma$. The sheaf~$\sP_\defpar$ becomes a sheaf of graded algebras  
by assigning degree $1$ to local sections~$x_i\in\sQ^\vee(U)$ 
as well as to~$z\in H^0(\O_C)$; thus its degree-zero piece is 
\[\sP_{\defpar,0}\cong \O_C\otimes\C\Gamma.
\]
\end{definition}

\begin{remark}\rm The sheaf of algebras~$\sA_\defpar$ is a relavitive version 
of the following non-commutative deformation of the skew group algebra, 
introduced by Crawley--Boevey and Holland in~\cite{cb-h}, depending on a 
deformation parameter $\lambda\in\lH\cong Z(\C\Gamma)$:
\[ A_\lambda= \C\langle x_1, x_2\rangle * \Gamma \big/\big\langle
\!\big\langle[x_1,x_2]+\lambda\big\rangle \!\big\rangle.
\]
The graded version is
\[ P_\lambda = \C\langle y_0, y_1, y_2\rangle * \Gamma \big/\big\langle \!\big\langle[y_0, y_i], [y_1, y_2]+\lambda y_0^2\big\rangle \!\big\rangle.\]

For~$\Gamma=\{1\}$, $\lambda$ is just a complex number; if $\lambda\neq 0$, 
$A_\lambda$ is isomorphic to the first Weyl algebra, whereas~$P_\lambda$ is 
a degenerate Sklyanin algebra deforming the algebra of functions on 
the commutative projective plane~$\PP^2$. As proved in~\cite{cb-h}, 
for general $\Gamma$ and $\lambda\in\lh\subset Z(\C\Gamma)$ the algebra
$A_\lambda$ is finite over its center
\[ Z A_\lambda \cong \C[\bar Y_\lambda]. 
\]
The latter is the coordinate ring of the affine variety $\bar Y_\lambda$ 
corresponding to the deformation parameter $\lambda\in\lh$, a deformation 
of the invariant ring $\C[x_1,x_2]^\Gamma\cong\C[\bar Y]$. For 
$\lambda\in \lH\setminus\lh$, $A_\lambda$ is ``genuinely'' 
non-commutative.
\end{remark}

By abuse of notation, we will refer to $\P_\defpar=\proj_C\sP_\defpar$ as the
non-commutative projective bundle corresponding to $\defpar\in \defspace$, with 
fibration $\pi_\dpar\colon\P_\defpar\rightarrow C$. 
Setting $z=0$, we have its divisor at infinity
\[i_\dpar\colon D_\defpar\hookrightarrow\P_\defpar.\] 
The divisor $D_\defpar$ has the structure of an ordinary (commutative) 
$\PP^1$-bundle
\[\pi_\dpar|_{D_\defpar}=\tau_\defpar\colon D_\defpar\rightarrow C\] 
equipped with a $\Gamma$-action on the fibres. Its complement 
$\A_\defpar=\P_\defpar\setminus D_\defpar = \spec_C \sA_\defpar$ is
a non-commutative affine bundle.

The sheaf $\sP_\defpar$ is a sheaf of regular graded algebras in the 
sense of~\cite{az}; sheaf theory on~$\P_\defpar$ works in 
complete analogy with the absolute case discussed in~\cite{bgk}. 
The category of coherent sheaves $\Coh(\P_\defpar)$ is by definition
the quotient of the category of sheaves of finitely generated 
graded right $\sP_\defpar$-modules by the subcategory of sheaves of 
torsion $\sP_\defpar$-modules; we will sometimes refer to objects in this
category as $\sP_\defpar$-modules. The trivial module, graded in
degree~$n$, defines the object $\O_{\P_\defpar}(n)\in\Coh(\P_\defpar)$; 
given a sheaf $\E$, its twists $\E(n)$ are obtained by shifting the grading. 
We have $\Ext$ groups as the derived functors of $\Hom$, and also 
functors $\shExt^i(-,\O_{\P_\defpar})$; the latter take values in the
category of left $\sP_\defpar$-modules (compare~\cite{bgk}). 

Pushforward 
\[\pi_{\defpar*}\colon \Coh(\P_\defpar)\rightarrow \Coh^\Gamma(C)\]
along the morphism $\pi_\dpar\colon\P_\defpar\rightarrow C$ is defined in 
the usual way, as the coherent $\Gamma$-sheaf on $C$ defined by sections over 
preimages of open sets of $C$, the section spaces being (right) $\C\Gamma$-modules; 
the action of $\Gamma$ on $C$ is taken to be trivial.
The higher pushforwards ${\rm R}^p\pi_{\defpar*}(-)$ are the derived functors 
of $\pi_{\defpar*}$. Given a $\P_\defpar$-module $\E$, we will 
also use the relative $\Hom$-functor
\[\Hom_C(\E, -)\colon \Coh(\P_\defpar)\rightarrow \Coh^\Gamma(C)\]
defined by homomorphisms on preimages of open sets in $C$, as well as its
derived functors 
\[\Ext^i_C(\E, -)\colon \Coh(\P_\defpar)\rightarrow \Coh^\Gamma(C).\]
We also have a pullback functor 
\[\pi_\dpar^*\colon \Coh^\Gamma(C)\rightarrow \Coh(\P_\defpar)\]
taking a sheaf of (right) $\C\Gamma$-modules $\F$ to the sheaf 
$\F\otimes_{\C\Gamma}\sP_\defpar$ of (right) $\sP_\defpar$-modules.
The pair $(\pi_\dpar^*,\pi_{\defpar*})$ forms an adjoint pair as 
in the commutative case. Similarly, for the inclusion 
$i_\dpar\colon D_\defpar\to\P_\defpar$, we have a pullback
(restriction) functor
\[i_\dpar^*\colon\coh(\P_\defpar)\rightarrow\coh^\Gamma(D_\defpar),\] 
defined by factoring modules of 
local sections by the ideal $\langle z\rangle$ (recall that $z$ is central),
as well as a pushforward 
\[i_{\defpar*}\colon\coh^\Gamma(D_\defpar)\rightarrow\coh(\P_\defpar),\] 
with $z$ acting on local sections by zero. 
There is also a restriction functor to the finite part $\A_\defpar$, defined
by factoring the ideal $\langle z-1\rangle$.

\begin{definition}\rm A $\pi_\dpar$-free sheaf on $\P_\defpar$ is an
object~$\E\in\Coh(\P_\defpar)$, which admits an embedding
\[  \E\hookrightarrow \pi_\dpar^*(\U)(n)
\]
for some $\U\in\Coh^\Gamma(C)$ and $n\in\Z$.
A framed $\pi_\dpar$-free sheaf $(\E, \phi)$ 
on~$(\P_\defpar, D_\defpar)$ is a $\pi_\dpar$-free sheaf~$\E$ on~$\P_\defpar$, 
together with a fixed isomorphism
\[ \phi\colon i_{\defpar}^*\E \stackrel{\sim}\longrightarrow \tau_\defpar^* \W,\]
on the divisor $D_\defpar$ at infinity, for some $\W\in\Coh^\Gamma(C)$.
\end{definition}

\begin{remark}\rm If $\pi\colon\P\rightarrow \{*\}$ is a (non-commutative) 
projective space over a point, the $\pi$-free sheaves are exactly the torsion free
ones (compare~\cite[Section 2]{bgk}). 
To see this, note that a $\pi$-free sheaf is certainly torsion free, 
since it embeds into a locally free sheaf. Conversely, a torsion free sheaf
embeds into some locally free sheaf, which in turn embeds into some~$\O^m_\P(n)$.
\end{remark}

\begin{lemma} If $\E$ is $\pi_\dpar$-free, then $L^ji_\dpar^*\E=0$ for $j>0$.
\end{lemma}
\begin{pf} As in the commutative case, the structure 
sheaf $i_{\defpar*}\O_{D_\defpar}$ has a resolution 
\[ 0 \to \O_{\P_\defpar}(-1) \stackrel{z}\to \O_{\P_\defpar}\to i_{\defpar*}\O_{D_\defpar}\to 0,\]
which implies that $L^ji_\dpar^*\E=0$ for $j>1$ for any $\E\in\Coh(\P_\defpar)$, 
and also that $L^1i_\dpar^*$ is left exact. If $\E$ is $\pi_\dpar$-free, 
applying the latter to an embedding $\E\hookrightarrow \pi_\dpar^*(\U)(n)$
gives the vanishing of $L^1$ also. 
\end{pf}


\section{Threefolds: the results}\label{sec!res} 

\subsection{Twisted quiver representations and quiver sheaves}
\label{sec!quiv}

Recall that, given a quiver with arrows $a\to b$ marked by 
objects $O_{ab}\in\sC$ of an abelian tensor category $\sC$, a representation 
of the marked quiver in $\sC$ consists of a set of objects $O_a$ of 
$\sC$ associated to 
nodes, and a set of morphisms $\phi_{ab}\in\Hom_\sC(O_a\otimes O_{ab}, O_b)$ 
associated to the arrows $a\to b$. Representations of a marked quiver in 
the category~$\Coh(X)$ of an algebraic variety~$X$ are also called
{\em quiver sheaves}~\cite{gk} on~$X$. 

\begin{figure}[ht]
\centering
\ifthenelse{\equal{\pics}{true}}
{\input{quiver1.pstex_t}}{}
\caption{The marked extended McKay quiver for $A_2$} 
\label{fig!quiv1}
\end{figure}
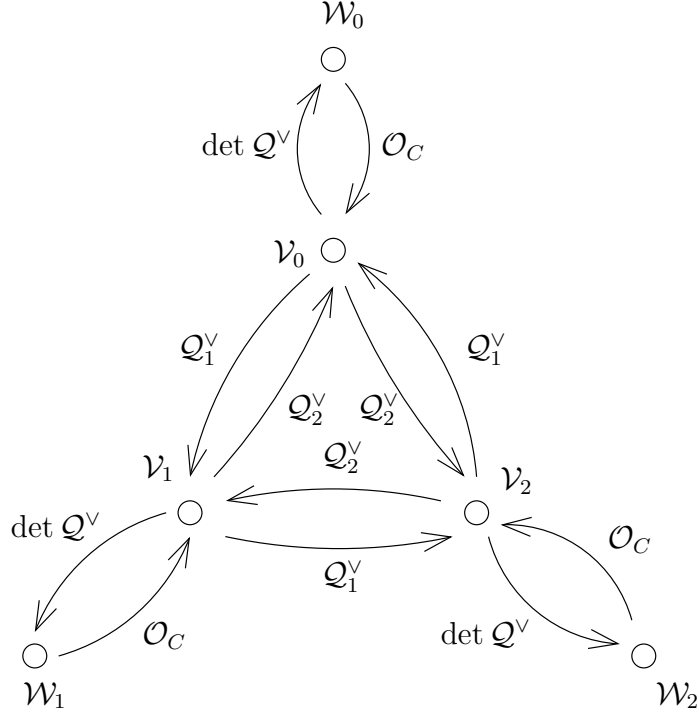

In the specific context of classifying holomorphic $D$-branes 
on the threefold~$X$ and its 
deformations, the following quiver marked in~$\Coh(C)$ will arise naturally. 
The quiver is the standard extended McKay quiver of~\cite{nak}, 
obtained from the original 
one by adding an extra leaf at each node with arrows in both directions.
Using the data of the vector bundle $\sQ$ on $C$, we mark this quiver 
in $\Coh(C)$ as follows.
\begin{itemize}
\item The marked $A_n$ quiver for $n>1$ is illustrated on 
Figure~\ref{fig!quiv1}; recall that in this case, there is a decomposition 
$\sQ=\sQ_1\oplus \sQ_2$ into a sum of line bundles, 
since the structure group of $\sQ$ reduces to the diagonal torus. 
\item The marked $A_1$ quiver consists of only two nodes $0$ and $1$ 
and two arrows $0\to 1, 1\to 0$ marked by the rank-two bundle $\sQ^\vee$, as well 
as leaves marked as in the higher $A_n$ case.  
\item For types $D$ and $E$, arrows between nodes are all marked by the line 
bundle $\sQ_0^\vee$, where $\sQ=\sQ_0^{\oplus 2}$; leaves are marked as before. 
\end{itemize}

\subsection{The main classification result}

\begin{theorem} Given $\defpar\in H^0(C, \sH)$, there is a 
1-to-1 correspondence between the following sets of data.
\label{thm!class}
\begin{enumerate} \renewcommand{\labelenumi}{\rm(\arabic{enumi})}
\item Isomorphism classes of framed $\pi_\dpar$-free sheaves $(\E, \phi)$
on $(\P_\defpar, D_\defpar)$.
\item Quintuples $(\V, \W, \B, \I, \J)$, where $\W, \V$ are coherent 
$\Gamma$-sheaves on $C$, and 
\[\begin{array}{rcl}
\B& \in&\Hom^\Gamma_C(\V\otimes \sQ^\vee, \V), \\ \\
\I &\in &\Hom^\Gamma_C(\W, \V),\\ \\ 
\J &\in&\Hom^\Gamma_C(\V\otimes \det\sQ^\vee, \W),\end{array}\] 
satisfying the following two conditions: 
\begin{enumerate}
\item the ADHM relation 
\[  \B\wedge \B + \I\circ\J + \defpar =0 \in \Hom^\Gamma_C(\V\otimes \det\sQ^\vee, \V),
\] 
where
\[H^0(C, Z(\C\Gamma)\otimes\det\sQ)\hookrightarrow \Hom^\Gamma_C(\V\otimes \det\sQ^\vee, \V)\]
is the natural embedding as the central subspace;
\item non-degeneracy: if $\V'\subset\V$ is a $\Gamma$-subsheaf such that 
$\B(\V'\otimes\sQ^\vee)\subset \V'$ and $\I\W\subset\V'$, then~$\V'=\V$. 
\end{enumerate} 
Sets of quintuples 
are identified under the action of invertible elements of $\Hom^\Gamma_C(\V,\V)$. 
\item Representations 
$(\{\V_a\}, \{\W_a\}, \{\B_{ab}\}, \{\I_a\}, \{\J_a\})$ 
in $\Coh(C)$ of the marked McKay-type qui\-ver introduced in~\ref{sec!quiv}, 
satisfying 
\begin{enumerate}
\item the ADHM relations
\[ \sum_b \epsilon_{ab}\B_{ba}\circ \B_{ab} + \I_a\circ \J_a + \defpar_a =0 \in \Hom_C(\V_a\otimes \det\sQ^\vee, \V_a)
\]
at each node~$a$, where $\epsilon_{ab}\in\{\pm 1\}$ is a standard assingment of 
signs to arrows with $\epsilon_{ab}=-\epsilon_{ba}$, and 
$\defpar_a=P_{\roo_a}(\defpar)$ is the projected section corresponding to 
the simple root $\roo_a$, and 
\item non-degeneracy: if $\{\V_a'\}$ is a $\B$-invariant set of subsheaves 
containing the images of $\I_a$'s, then $\V_a'=\V_a$ at all nodes. 
\end{enumerate}
Two representations
are identified under invertible elements of $\prod_a\Hom_C(\V_a,\V_a)$.
\end{enumerate}

\noindent If $\defpar\in\defspace=H^0(C, \sh)$ is a deformation parameter of
the threefold $X=X_0$, then the same data also parametrizes
\begin{enumerate} \renewcommand{\labelenumi}{\rm(\arabic{enumi})}
\addtocounter{enumi}{3}
\item certain objects in $\D(\Coh X_\defpar)$, the derived category 
of coherent sheaves on $X_\defpar$. 
\end{enumerate}
\end{theorem} 
\begin{pf}
The equivalence ${\rm (1)}\iff {\rm (2)}$ follows from a version of 
the relative Beilinson resolution for the non-commutative 
projective bundle $\P_\defpar\rightarrow C$; details are given in 
Section~\ref{sec!koszul}. 
McKay's definition of the quiver describing the representation theory 
of $\Gamma$ implies ${\rm (2)}\iff {\rm (3)}$ in the standard way.
Finally the mapping ${\rm (1)}\implies {\rm (4)}$ in the geometric 
case~$\defpar\in\defspace=H^0(C, \sh)$ is given by a derived 
equivalence to be discussed in Section~\ref{sec!derived}.
\end{pf}

\begin{remark}\rm As $X=X_0$ and its deformations $X_\defpar$ for 
$\defpar\in \defspace$ are not projective, one needs to rigidify 
before holomorphic $D$-branes, in other words objects in $\D^b(X_\defpar)$, 
have a sensible moduli space. 
For the central fibre $X=X_0$, a crepant resolution of the singular threefold
$\sQ/\Gamma$, one has a derived equivalence~\cite{bkr}
\[ \D(X_0)\cong \D^\Gamma(\sQ)
\]
between the derived categories of coherent sheaves on~$X_0$ and that of
$\Gamma$-equivariant sheaves on the total space of the bundle~$\sQ\to C$. 
One can easily rigidify
on the latter by considering~$\Gamma$-sheaves on the 
projective bundle $\P_0=\PP(\sQ\oplus\O_C)\rightarrow C$, framed on the divisor 
at infinity $D_0 =\PP(\sQ)\hookrightarrow\P_0$.
Theorem~\ref{thm!class} is the appropriate generalization of this approach 
which also works for deformations: for the analogous problem on $X_\defpar$, we 
consider framed sheaves on the non-commutative projective 
bundle~$\P_\defpar\rightarrow C$.

In the surface case, this approach was used earlier in~\cite{bgk}. 
To quote the result, let $\lambda\in Z(\C\Gamma)$. 
Then for $\Gamma$-modules $V,W$, Nakajima's non-singular quiver 
variety $\M_{V, W, \lambda}$ parametrizes torsion free
sheaves on the non-commutative space $\PP^2_\lambda=\proj P_\lambda$, 
framed on the commutative $\Gamma$-line at $\infty$. This statement 
generalizes earlier work of \cite{don, kn, nak1, n-s, kko}
and others. The origin of all such results is of course the ADHM 
classification~\cite{adhm} of finite-action $\SU(\dim(W))$-instantons 
on~$\R^4$ of charge $\dim(V)$.
\end{remark}

\subsection{Some holomorphic $D$-branes on $ADE$ fibrations over $\AA^1$}

If $C\cong\AA^1$, Theorem~\ref{thm!class} can in some cases be re-written in 
terms of classical quiver representations: representations of a quiver in 
vector spaces. This will give an interpretation of an assertion
of~\cite{ckv, katz, zhu}.

Recall that for $C\cong\AA^1$, a deformation parameter $\defpar\in\defspace$ 
of the central fibre $X_0=\AA^1\times Y$ can be specified by a 
set of polynomials $\{\defpoly_a\in \C[t] : a\in\De\}$ indexed by the vertices of 
the affine quiver, subject to $\sum_a \delta_a\cdot \defpoly_a = 0$. The 
following definition is due to Cachazo--Katz--Vafa~\cite{ckv, katz}. 

\begin{definition}\rm The {affine \susy\ $ADE$ quiver} is the McKay quiver 
extended by a loop $a\to a$ at each vertex. For a (finite-dimensional)
representation $(\{V_a\}, \{B_{ab}\}, \{\Psi_a\})$ of this quiver, 
the ADHM-type relations are 
\begin{equation}
\label{rel1}
 \sum_{b} \epsilon_{ab} B_{ba} B_{ab} + \defpoly_a(\Psi_a)=0\in \Hom(V_a, V_a)
\end{equation}
at each vertex $a\in\De$ of the quiver, where $\defpoly_a(\Psi_a)$ is to be 
interpreted as the evaluation of a polynomial on an endomorphism of $V_a$, 
as well as
\begin{equation}
\label{rel2}
\Psi_a B_{ba} = B_{ba}\Psi_b\in\Hom(V_a, V_b)
\end{equation}
along each arrow $a\to b$ of the quiver~$\De$.
\end{definition}

Consider quadruples $(\{V_a\}, \{B_{ab}\}, \{\Psi_a\}, \spvec)$, 
where $(\{V_a\}, \{B_{ab}\}, \{\Psi_a\})$ is a representation of the 
affine \susy\ $ADE$ quiver satisfying the ADHM-type relations, 
and $\spvec\in V_0$ is a fixed vector in the vector space attached to the 
affine node. Call a quadruple non-degenerate if
there is no $(B, \Psi)$-invariant collection of subspaces $\{V_a'\subset V_a\}$
with~$\spvec\in V_0'$. 

\begin{theorem} Equivalence classes of non-degenerate quadruples 
$(\{V_a\}, \{B_{ab}\}, \{\Psi_a\}, \spvec)$ satisfying the ADHM relations, 
identified under the action of $\prod_a \GL(V_a)$, parametrize certain objects 
in $\D(\coh X_\defpar)$, holomorphic $D$-branes on the threefold~$X_\defpar$. 
\label{thm!A1}
\end{theorem}

\begin{pf} Quiver sheaf data on $C$ parametrize certain branes on~$X_\defpar$ 
by Theorem~\ref{thm!class}. The correspondence between representations 
of the \susy\ $ADE$ quiver and a special class of quiver sheaf data will 
be discussed in Section~\ref{sec!affpf}. 
\end{pf}

\begin{remark}\rm As explained in~\cite{ckv}, the quiver 
relations~(\ref{rel1})-(\ref{rel2}) come from the natural superpotential 
of the quiver gauge theory on $\De$, involving adjoint 
fields $\Psi_a$ as well as bifundamental fields $B_{ab}$. 
\end{remark}

\begin{remark}\rm Let the finite \susy\ $ADE$ quiver be obtained from the affine 
one by deleting the affine node. Representations of the finite \susy\ $ADE$ quiver, 
satisfying the ADHM-type relations~(\ref{rel1})-(\ref{rel2}), 
parametrize holomorphic $D$-branes supported on exceptional fibres of
$f_\dpar\colon X_\defpar\rightarrow \bar X_\defpar$. This follows from the
statement that the vanishing of the affine component of $\V$ forces 
all other $\V_a$ be supported on points $P\in C$ at which some projected section
$p_\roo(\defpar)$ vanishes for some positive root $\roo\in\posro$, 
in other words on points of the base curve over which 
the surface fiber $\pi_\dpar^{-1}(P)$ contains exceptional curves.
Observing that the section $\defpar\in H^0(C, \Z(\C\Gamma)\otimes\det\sQ)$ 
is central in $\Hom^\Gamma_C(\V\otimes \det\sQ^\vee, \V)$, 
so commutes with all components of $\B$, the latter statement is essentially 
proved in~\cite[4.1--4.2]{ckv}.
This establishes a direct link to~\cite{katz, zhu}, according to 
which (in the generic case) irreducible representations of the finite \susy\ quiver
with the given relations parametrize exceptional components of 
the contraction $f_\dpar\colon X_\defpar\rightarrow \bar X_\defpar$.
\end{remark}

\section{Proofs}\label{sec!pfs} 
\subsection{The Beilinson argument} 
\label{sec!koszul}

The aim of this section is to prove of the equivalence 
${\rm (1)}\iff {\rm (2)}$ of the classification result Theorem~\ref{thm!class} 
via an analysis of framed $\pi_\dpar$-free sheaves on $\P_\defpar$. 

Given $\defpar\in H^0(C, \sH)$, recall the sheaf of algebras~$\sP_\defpar$ 
on the curve~$C$, and the associated non-commutative bundle
$\pi_\dpar\colon\P_\defpar\rightarrow C$. 
Define $\sP_\defpar$-modules~$\T_i$ by 
\begin{equation}
\begin{array}{c} 
\T_0 = \O_{\P_\defpar},\\ \\
0 \longrightarrow \O_{\P_{\defpar}} \longrightarrow \pi_\dpar^*(\sQ\oplus\O_C)(1) \longrightarrow \T_1 \longrightarrow 0, \\ \\
\T_2 =  \pi_\dpar^* (\det\sQ)(3).
\label{seq!defN1}
\end{array}
\end{equation}

\begin{proposition} 
A $\pi_\dpar$-free sheaf $\E$ on $\P_\defpar$, framed on the divisor $D_\defpar$, 
is the cohomology of a monad
\[ \pi_\dpar^*\Ext^1_C\left(\T_2(-1), \E\right)(-1)\rightarrow \pi_\dpar^*\Ext^1_C\left(\T_1, \E\right)\rightarrow\pi_\dpar^*\Ext^1_C\left(\T_0(1), \E\right)(1)
\]
 of $\sP_\defpar$-modules.
\label{prop!Emon}
\end{proposition}
\begin{pf} 
Given a $\sP_\defpar$-module $\F$, a Koszul duality argument, in an analogous way 
to the absolute case in~\cite[Section 7]{bgk} following~\cite[Thm 2.6.1]{bgs}, 
leads to a Beilinson-type spectral sequence with $E_1$ term
\[ E_1^{p,q} = \pi_\dpar^*\Ext^q_C\left(\T_{-p}(p), \F\right)(p),
\] 
nonzero only for $-2\leq p\leq 0, \ 0\leq q \leq 2$, converging to~$\F$ 
in the limit. The vanishing results 
\[ \Ext^q_C\left(\T_{-p}(p), \E(-1)\right) =0 \mbox{ for } q=0,2,\ p=-1, -2
\]
which follow from the existence of the framing of $\E$ on the divisor $D_\defpar$
(compare~\cite[Lemma 6.2]{kko},~\cite[Lemma 4.2.12]{bgk}),
reduce the spectral sequence for $\F=\E(-1)$ to the monad given in the statement. 
Details are left to the reader.
\end{pf}

We also record an auxiliary result. 

\begin{lemma}
\label{lem!natiso}
There are natural isomorphisms
\[\Hom_{\P_\defpar}\!\!\left(\pi_\dpar^*\det\sQ^\vee, \pi_\dpar^*\sQ^\vee(1)\right) \cong \Hom_{\P_\defpar}\!\!\left(\pi_\dpar^*\sQ^\vee, \O_{\P_\defpar}(1)\right) \cong \Hom^\Gamma_C\!\left(\sQ^\vee,\sQ^\vee\oplus \O_C\right).
\]
\end{lemma}

\begin{pf} The first isomorphism follows from Lemma~\ref{lem!rk2} below. The 
second one follows from adjunction for the pair~$(\pi_\dpar^*,\pi_{\defpar*})$, 
together with 
\[\pi_{\defpar*}\O_{\P_\defpar}(1)\cong \sP_{\defpar,1}\cong (\sQ^\vee\oplus \O_C)*\C\Gamma\in\Coh^\Gamma(C),\]
an identity well known from the commutative context. 
\end{pf}

\begin{lemma} Let $\sQ$ be a rank-two bundle on a (commutative) space. 
Then there is a natural isomorphism
\[ \sQ\otimes\det\sQ^\vee \cong \sQ^\vee.
\] 
\label{lem!rk2}
\end{lemma}
\begin{pf} The embedding 
$\iota\colon\det\sQ^\vee\rightarrow \sQ^\vee\otimes\sQ^\vee$ 
induces a natural map
\[\Hom\!\left(\det\sQ^\vee, \det\sQ^\vee\right)\rightarrow\Hom\!\left(\det\sQ^\vee, (\sQ^\vee)^{\otimes 2}\right)\cong \Hom\!\left(\sQ\otimes\det\sQ^\vee, \sQ^\vee\right). 
\]
The image of the identity of the first $\Hom$-group gives a natural morphism 
as in the statement, which can be checked on a local basis to 
be an isomorphism.
\end{pf}

Now return to the context of the classification result Theorem~\ref{thm!class}, 
and consider a quintuple $(\V, \W, \B, \I, \J)$ as in Theorem~\ref{thm!class}(2); 
recall that \[\W, \V\in\Coh^\Gamma(C),\] and 
\[\begin{array}{rcl}
\B& \in&\Hom^\Gamma_C(\V\otimes \sQ^\vee, \V), \\ \\
\I &\in &\Hom^\Gamma_C(\W, \V),\\ \\ 
\J &\in&\Hom^\Gamma_C(\V\otimes \det\sQ^\vee, \W).\end{array}\] 
Let
\[ c\in\Hom_{\P_\defpar}\!\!\left(\pi_\dpar^*\det\sQ^\vee, \pi_\dpar^*\sQ^\vee(1)\right), \ \ \ d\in \Hom_{\P_\defpar}\!\!\left(\pi_\dpar^*\sQ^\vee, \O_{\P_\defpar}(1)\right)
\]
denote the images, under the isomorphisms of~Lemma~\ref{lem!natiso}, of the 
canonical element 
\[{\rm Id}\in\Hom^\Gamma_C( \sQ^\vee,\sQ^\vee)\subset\Hom^\Gamma_C\!\left(\sQ^\vee,\sQ^\vee\oplus \O_C\right)\]
Note also that we have a fixed section 
\[z\in\Hom_{\P_\defpar}(\O_{\P_\defpar}, \O_{\P_\defpar}(1)).\] 
Define 
\[
a = \left(\begin{array}{c} \pi_\dpar^*\left(\B\circ({\rm Id}_\V\otimes\iota)\right)\otimes z - \pi_\dpar^*({\rm Id_\V})\otimes c\,(-1) \\ \pi_\dpar^*(\J)\otimes z \end{array} \right) \colon \pi_\dpar^*(\V\otimes\det\sQ^\vee)(-1) \longrightarrow \pi_\dpar^*(\V\otimes\sQ^\vee\oplus \W)
\]
where $\iota\colon \det\sQ^\vee\to(\sQ^\vee)^{\otimes 2}$ is the natural map.
Define similarly
\[
b= \left(\begin{array}{cc}\pi_\dpar^*(\B)\otimes z + \pi_\dpar^*({\rm Id_\V})\otimes d & \pi_\dpar^*(\I)\otimes z \end{array}\right) :  \pi_\dpar^*(\V\otimes\sQ^\vee\oplus \W)\to \pi_\dpar^*(\V)(1),
\]
to obtain the chain of morphisms
\begin{equation}
\pi_\dpar^*(\V\otimes\det\sQ^\vee)(-1)\stackrel{a}\longrightarrow \pi_\dpar^*(\V\otimes\sQ^\vee\oplus \W)\stackrel{b}\longrightarrow \pi_\dpar^*(\V)(1).
\label{seq!monad}
\end{equation}

The following result completes the proof of the equivalence 
${\rm (1)}\iff {\rm (2)}$ of the classification result Theorem~\ref{thm!class}.

\begin{proposition} If the quintuple satisfies the ADHM relation, then
{\rm (\ref{seq!monad})}~is a complex of $\sP_\defpar$-modules.
Furthermore, it is a monad defining a framed $\pi_\dpar$-free sheaf $\E$ if
and only if the quintuple  $(\V, \W, \B, \I, \J)$ is non-degenerate.
Conversely, every $\pi_\dpar$-free $\sP_\defpar$-module~$\E$, framed 
on~$D_\defpar$, arises from this construction.
\end{proposition}
\begin{pf} The standard direct computation shows that $b\circ a=0$ is equivalent
to the ADHM relation. The proof of the equivalence of the monad property and 
non-degeneracy is analogous to the absolute case~\cite[Section 4.1]{bgk}. 
For the converse, given a framed sheaf $(\E, \phi)$, 
let $\V=\Ext^1_C(\O_{\P_\defpar}(1), \E)$. Then by Proposition~\ref{prop!Emon}, 
$\E$ is the middle cohomology of the monad
\[ \pi_\dpar^*(\V\otimes\det\sQ^\vee)(-1)\rightarrow \pi_\dpar^*\Ext^1_C\left(\T_1, \E\right)\rightarrow\pi_\dpar^*\V(1).
\]
The usual arguments~\cite[Theorem 6.7]{kko} show that, since $\E$ is framed
on~$D_\defpar$, this monad is isomorphic to a monad of the form~(\ref{seq!monad}) 
for some quintuple $(\V, \W, \B, \I, \J)$.
\end{pf}

\subsection{A derived equivalence}\label{sec!derived}

In this section we complete the proof of Theorem~\ref{thm!class} by establishing
the missing link ${\rm (1)}\implies {\rm (4)}$. 

\begin{proposition} Let $\defpar\in \defspace$ be a deformation 
parameter of the central fibre $X=X_0$.
There is a distinguished equivalence of triangulated categories 
\[ \D(\coh X_\defpar) \cong \D(\mod\sA_\defpar ),
\] 
where $\mod\sA_\defpar$ is the category of sheaves of 
finitely generated right $\sA_\defpar$-modules, and $\D(-)$ denotes the bounded 
derived category on both sides.
\end{proposition}
\begin{pf} This assertion is a fibered version of the analogous two-dimensional 
equivalence proved in~\cite{gs}, and the proof carries over verbatim.
A deformation argument starting from 
the central fibre $X=X_0$ shows that a certain specific component $M_\defpar$ 
of a fine moduli space of torsion sheaves on $\A_\defpar$ maps by a 
semi-small birational map to the singular variety $\bar X_\defpar$. 
By~\cite{vdb}, generalizing an argument of~\cite{bkr},
this implies that $M_\defpar$ is a crepant resolution of $\bar X_\defpar$, 
and one has a derived equivalence
\[ \D(\coh M_\defpar) \cong \D(\mod\sA_\defpar)
\]
defined by the universal sheaf.  
But since $X_\defpar$ is the unique crepant resolution of $\bar X_\defpar$, 
necessarily $M_\defpar\cong X_\defpar$ and the proposition follows. 
Details are left to the reader.
\end{pf}

This equivalence gives the mapping ${\rm (1)}\implies {\rm (4)}$ 
of Theorem~\ref{thm!class} from framed $\pi_\dpar$-free sheaves on $\P_\defpar$
to objects in $\D(\coh X_\defpar)$. Indeed, a right $\sP_\defpar$-module 
can be restricted to the affine part $\A_\defpar$ to give a right
$\sA_\defpar$-module, and then mapped using the derived equivalence to an object
in $\D(\coh X_\defpar)$, in other words a holomorphic $D$-brane on $X_s$.

\subsection{Fibrations over the affine line}
\label{sec!affpf}

In this section, we take a fibration $X_\defpar\rightarrow C\cong\AA^1$ and 
discuss the proof of Theorem~\ref{thm!A1}. From Theorem~\ref{thm!class}, 
we know that certain holomorphic $D$-branes on $X_\defpar$ are classified by 
non-degenerate quintuples $(\V,\W, \B, i, j)$ satisfying the ADHM equation.
Consider the subclass of representations in $\Coh(\AA^1)$ with the simplest 
possible framing
$\W\cong\O_{\AA^1}$ and $\V$ a torsion $\Gamma$-sheaf on $\AA^1$. It follows
that $\J=0$ and $\I\in H^0(\AA^1, \V^\Gamma)$.
Decompose $\V$ and the map $B$ into $\Gamma$-components to obtain
torsion sheaves $\V_a$ and sheaf homomorphisms 
$\B_{ab}\colon \V_a\rightarrow \V_b$
indexed by nodes and edges of the McKay quiver.

Set $V_a=H^0(\AA^1, \V_a)$, and let $B_{ab}=H^0(\B_{ab})\colon V_a\rightarrow V_b$ 
to be the map on global sections induced by $\B_{ab}$. Let 
$\spvec\in V_0$ be the section corresponding to $\I$.
Let also $\Psi_a\colon V_a\rightarrow V_a$ be the map induced by 
multiplication by the section $t\in H^0(\AA^1, \O_{\A_1}) \cong \C[t]$. 
Theorem~\ref{thm!A1} follows from Theorem~\ref{thm!class}, together with

\begin{proposition} The map
\[ (\V, \O_C, \B, 0, 0) \mapsto (\{V_a\}, \{B_{ab}\}, \{\Psi_a\}, \spvec\in V_0)
\] 
sets up a one-to-one correspondence from this restricted set
of quiver ADHM data to representations of the affine \susy\ ADE quiver 
satisfying the relations~(\ref{rel1})-(\ref{rel2}).
\end{proposition} 
\begin{pf} Given $(\V, \B)$, the edge relations~(\ref{rel2}) 
$\Psi_a B_{ba} = B_{ba}\Psi_b$ 
for the data $(\{V_a\}, \{B_{ab}\}, \{\Psi_a\})$ hold by definition. 
Further, the ADHM equation for $(\V, \B)$ is
\[ \B\wedge \B + \defpar = 0 \in \Hom(\V, \V\otimes\det\sQ),\]
which in $\Gamma$-components says that 
\[ \sum_{b}\epsilon_{ab}\B_{ba}\circ\B_{ab}+\defpar_a=0\in \Hom(\V_a, \V_a).
\]
Replacing $\defpar_a$ by the polynomial $\defpoly_a$, and remembering
that the effect of $t\in H^0(\O_{\AA^1})$ on $H^0(\V)$ is exactly $\Psi_a$, 
for global sections we obtain
\[\sum_{b} \epsilon_{ab} B_{ba}\circ B_{ab}+\defpoly_a(\Psi_a)=0\in \Hom(V_a, V_a)
\]
which is exactly relation~(\ref{rel1}) for the node~$a$. 

Conversely, given a representation $(\{V_a\}, \{B_{ab}\}, \{\Psi_a\},\spvec\in V_0)$
of the \susy\ ADE quiver, define torsion sheaves attached to the nodes by
\[\begin{CD}\V_a = \coker\Big(V_a\otimes\O_{\AA^1}@>{1\otimes t - \Psi_a\otimes 1}>> V_a\otimes \O_{\AA^1}\Big).\end{CD}\]
Using Lemma~\ref{linalg} below, for adjacent nodes $a,b$ we have a diagram
\[\begin{CD} 0 & @>>> & {\ \ V_a\otimes\O_{\AA^1}} & @>{1\otimes t - \Psi_a\otimes 1}>> & {\ \ V_a\otimes\O_{\AA^1}} & @>>> \V_a @>>> 0 \\
& & & & @V{B_{ab}\otimes 1}VV && @V{B_{ab}\otimes 1}VV \\
0 & @>>> & {\ \ V_b\otimes\O_{\AA^1}} & @>{1\otimes t - \Psi_b\otimes 1}>> & {\ \ V_b\otimes\O_{\AA^1}} & @>>> \V_b @>>> 0
\end{CD}
\]
which, by commutativity $\Psi_a B_{ba} = B_{ba}\Psi_b$, induces a map
$\B_{ab}\colon \V_a\rightarrow \V_b$. The converse of the above argument 
shows that the ADHM relation follows from the relations~(\ref{rel1}). 
By Lemma~\ref{linalg}, the two constructions are inverses to each other. 
\end{pf}

The proof used the elementary

\begin{lemma} Given a torsion sheaf $\V$ on $\AA^1=\spec \C[t]$, 
let $V=H^0(\AA^1, \V)$ and let $\Psi\colon V\rightarrow V$ be the map 
given by multiplication by $t\in H^0(\O_{\AA^1})$. Then the sequence of sheaves 
\[\begin{CD} 0  @>>> V\otimes\O_{\AA^1} @>{1\otimes t - \Psi\otimes 1}>>  V\otimes\O_{\AA^1}  @>c>> \V @>>> 0\end{CD}\]
is exact on~$\AA^1$, where $c\colon H^0(\V)\otimes \O_{\AA^1} \rightarrow \V$ is 
the canonical map. Conversely, given a vector space with an endomorphism
$(V, \Psi)$, the exact sequence defines a torsion sheaf $\V$ on $\AA^1$, 
and the two constructions are mutual inverses. 
\label{linalg}\end{lemma} 

\begin{remark}\rm In this Lemma, $\V\cong\O_Z$ is a structure sheaf of 
a $0$-dimensional subscheme $Z\subset\AA^1$ if and only if $\Psi$ is a
regular endomorphism. Their moduli space is  
\[{\rm Mat}(n,\C)/\!/\GL(n,\C)\cong\{\mbox{regular endomorphisms}\}/\GL(n,\C)\cong \AA^n\cong(\AA^1)^{[n]},\] 
where the map is given by taking the coefficients of the characteristic polynomial 
of $\Psi$, which is also the equation of the corresponding subscheme.
\end{remark}
\subsection*{Acknowledgements} Thanks to Sheldon Katz, Eduard Looijenga, 
Tom Nevins and Tony Pantev for helpful remarks and correspondence. 
Special thanks to Ian Grojnowski for many conversations on subjects 
related to this paper. Support by a European Union Marie Curie 
Individual Fellowship and by OTKA grant $\#046878$ is also gratefully
acknowledged.

\noindent {\small {\sc Department of Mathematics, Utrecht University}

\noindent Current address:

\noindent {\sc Mathematical Institute, University of Oxford}

\noindent E-mail address: \tt szendroi@maths.ox.ac.uk}

\end{document}

%% file: quiver1.pstex_t
\begin{picture}(0,0)%
\includegraphics{quiver1.pstex}%
\end{picture}%
\setlength{\unitlength}{3947sp}%
\begingroup\makeatletter\ifx\SetFigFont\undefined%
\gdef\SetFigFont#1#2#3#4#5{%
  \reset@font\fontsize{#1}{#2pt}%
  \fontfamily{#3}\fontseries{#4}\fontshape{#5}%
  \selectfont}%
\fi\endgroup%
\begin{picture}(4058,4473)(376,-3715)
\put(4426,-3661){\makebox(0,0)[lb]{\smash{\SetFigFont{12}{14.4}{\familydefault}{\mddefault}{\updefault}{\color[rgb]{0,0,0}$\W_2$}%
}}}
\put(3451,-2311){\makebox(0,0)[lb]{\smash{\SetFigFont{12}{14.4}{\familydefault}{\mddefault}{\updefault}{\color[rgb]{0,0,0}$\V_2$}%
}}}
\put(1201,-2236){\makebox(0,0)[lb]{\smash{\SetFigFont{12}{14.4}{\familydefault}{\mddefault}{\updefault}{\color[rgb]{0,0,0}$\V_1$}%
}}}
\put(2026,-886){\makebox(0,0)[lb]{\smash{\SetFigFont{12}{14.4}{\familydefault}{\mddefault}{\updefault}{\color[rgb]{0,0,0}$\V_0$}%
}}}
\put(2326,614){\makebox(0,0)[lb]{\smash{\SetFigFont{12}{14.4}{\familydefault}{\mddefault}{\updefault}{\color[rgb]{0,0,0}$\W_0$}%
}}}
\put(3226,-1486){\makebox(0,0)[lb]{\smash{\SetFigFont{12}{14.4}{\familydefault}{\mddefault}{\updefault}{\color[rgb]{0,0,0}$\sQ_1^\vee$}%
}}}
\put(2101,-1861){\makebox(0,0)[lb]{\smash{\SetFigFont{12}{14.4}{\familydefault}{\mddefault}{\updefault}{\color[rgb]{0,0,0}$\sQ_2^\vee$}%
}}}
\put(2326,-2911){\makebox(0,0)[lb]{\smash{\SetFigFont{12}{14.4}{\familydefault}{\mddefault}{\updefault}{\color[rgb]{0,0,0}$\sQ_1^\vee$}%
}}}
\put(2551,-1861){\makebox(0,0)[lb]{\smash{\SetFigFont{12}{14.4}{\familydefault}{\mddefault}{\updefault}{\color[rgb]{0,0,0}$\sQ_2^\vee$}%
}}}
\put(2701,-211){\makebox(0,0)[lb]{\smash{\SetFigFont{12}{14.4}{\familydefault}{\mddefault}{\updefault}{\color[rgb]{0,0,0}$\O_C$}%
}}}
\put(1201,-3286){\makebox(0,0)[lb]{\smash{\SetFigFont{12}{14.4}{\familydefault}{\mddefault}{\updefault}{\color[rgb]{0,0,0}$\O_C$}%
}}}
\put(3076,-3286){\makebox(0,0)[lb]{\smash{\SetFigFont{12}{14.4}{\familydefault}{\mddefault}{\updefault}{\color[rgb]{0,0,0}$\det\sQ^\vee$}%
}}}
\put(4126,-2686){\makebox(0,0)[lb]{\smash{\SetFigFont{12}{14.4}{\familydefault}{\mddefault}{\updefault}{\color[rgb]{0,0,0}$\O_C$}%
}}}
\put(1426,-1486){\makebox(0,0)[lb]{\smash{\SetFigFont{12}{14.4}{\familydefault}{\mddefault}{\updefault}{\color[rgb]{0,0,0}$\sQ_1^\vee$}%
}}}
\put(376,-2611){\makebox(0,0)[lb]{\smash{\SetFigFont{12}{14.4}{\familydefault}{\mddefault}{\updefault}{\color[rgb]{0,0,0}$\det\sQ^\vee$}%
}}}
\put(1576,-211){\makebox(0,0)[lb]{\smash{\SetFigFont{12}{14.4}{\familydefault}{\mddefault}{\updefault}{\color[rgb]{0,0,0}$\det\sQ^\vee$}%
}}}
\put(2326,-2161){\makebox(0,0)[lb]{\smash{\SetFigFont{12}{14.4}{\familydefault}{\mddefault}{\updefault}{\color[rgb]{0,0,0}$\sQ_2^\vee$}%
}}}
\put(451,-3661){\makebox(0,0)[lb]{\smash{\SetFigFont{12}{14.4}{\familydefault}{\mddefault}{\updefault}{\color[rgb]{0,0,0}$\W_1$}%
}}}
\end{picture}

%% file: quiver2.pstex_t
\begin{picture}(0,0)%
\includegraphics{quiver2.pstex}%
\end{picture}%
\setlength{\unitlength}{3947sp}%
\begingroup\makeatletter\ifx\SetFigFont\undefined%
\gdef\SetFigFont#1#2#3#4#5{%
  \reset@font\fontsize{#1}{#2pt}%
  \fontfamily{#3}\fontseries{#4}\fontshape{#5}%
  \selectfont}%
\fi\endgroup%
\begin{picture}(3075,3048)(751,-3040)
\put(2101,-1861){\makebox(0,0)[lb]{\smash{\SetFigFont{12}{14.4}{\familydefault}{\mddefault}{\updefault}{\color[rgb]{0,0,0}$B_{10}$}%
}}}
\put(2551,-1861){\makebox(0,0)[lb]{\smash{\SetFigFont{12}{14.4}{\familydefault}{\mddefault}{\updefault}{\color[rgb]{0,0,0}$B_{02}$}%
}}}
\put(2326,-2236){\makebox(0,0)[lb]{\smash{\SetFigFont{12}{14.4}{\familydefault}{\mddefault}{\updefault}{\color[rgb]{0,0,0}$B_{21}$}%
}}}
\put(2326,-2986){\makebox(0,0)[lb]{\smash{\SetFigFont{12}{14.4}{\familydefault}{\mddefault}{\updefault}{\color[rgb]{0,0,0}$B_{12}$}%
}}}
\put(1426,-1486){\makebox(0,0)[lb]{\smash{\SetFigFont{12}{14.4}{\familydefault}{\mddefault}{\updefault}{\color[rgb]{0,0,0}$B_{01}$}%
}}}
\put(3076,-1411){\makebox(0,0)[lb]{\smash{\SetFigFont{12}{14.4}{\familydefault}{\mddefault}{\updefault}{\color[rgb]{0,0,0}$B_{20}$}%
}}}
\put(2026,-886){\makebox(0,0)[lb]{\smash{\SetFigFont{12}{14.4}{\familydefault}{\mddefault}{\updefault}{\color[rgb]{0,0,0}$V_0$}%
}}}
\put(2326,-136){\makebox(0,0)[lb]{\smash{\SetFigFont{12}{14.4}{\familydefault}{\mddefault}{\updefault}{\color[rgb]{0,0,0}$\Psi_0$}%
}}}
\put(751,-2986){\makebox(0,0)[lb]{\smash{\SetFigFont{12}{14.4}{\familydefault}{\mddefault}{\updefault}{\color[rgb]{0,0,0}$\Psi_1$}%
}}}
\put(3376,-2311){\makebox(0,0)[lb]{\smash{\SetFigFont{12}{14.4}{\familydefault}{\mddefault}{\updefault}{\color[rgb]{0,0,0}$V_2$}%
}}}
\put(1201,-2311){\makebox(0,0)[lb]{\smash{\SetFigFont{12}{14.4}{\familydefault}{\mddefault}{\updefault}{\color[rgb]{0,0,0}$V_1$}%
}}}
\put(3826,-2986){\makebox(0,0)[lb]{\smash{\SetFigFont{12}{14.4}{\familydefault}{\mddefault}{\updefault}{\color[rgb]{0,0,0}$\Psi_2$}%
}}}
\end{picture}